\documentclass[12pt]{amsart}
\usepackage{amsopn}
\usepackage{amssymb, amscd, verbatim}
\usepackage{graphics, epsfig}

\swapnumbers
\theoremstyle{plain}
\newtheorem{theorem}{Theorem}[section]
\newtheorem{proposition}[theorem]{Proposition}
\newtheorem{corollary}[theorem]{Corollary}
\newtheorem{lemma}[theorem]{Lemma}
\newtheorem*{teor}{Theorem}

\theoremstyle{definition}
\newtheorem{definition}[theorem]{Definition}

\newtheorem{remark}[theorem]{Remark}
\newtheorem{remarks}[theorem]{Remarks}
\newtheorem{nada}[theorem]{}
\newtheorem{notandrem}[theorem]{Notation and Remarks}

\theoremstyle{remark}
\newtheorem*{weakremark}{Remark}
\newtheorem*{weakremarks}{Remarks}

\def\R{\mathbf R}
\def\Z{\mathbf Z}
\def\C{\mathbf C}
\def\spec{\text{spec}}
\def\Lap{\Delta}
\def\bs{{\backslash}}

\def\O{{\mathcal O}}

\def\rt{{\R^2}}

\def\spm{{\spec_m}}

\def\A{{\mathcal A}}

\def\tm{{\langle\tau\rangle\bs M}}
\def\sm{{\langle\sigma\rangle\bs M}}
\def\ml{{m_\lambda}}
\def\hl{{H_\lambda}}
\def\hpl{{H^+_\lambda}}
\def\hml{{H^-_\lambda}}
\def\P{{{\bold P}}}
\def\tb{{\tau_P}}
\def\th{{\tau_H}}
\def\tvo{{\tau_{V1}}}
\def\tvt{{\tau_{V2}}}
\def\tone{{\tau_1}}
\def\ttwo{{\tau_2}}
\def\tthree{{\tau_3}}
\def\tfour{{\tau_4}}

\title[What the middle degree Hodge spectrum doesn't reveal]
{Boundary volume and length spectra of Riemannian manifolds:
What the middle degree Hodge spectrum doesn't reveal.}

\author[C.S. Gordon \and J.P. Rossetti]{Carolyn S.\ Gordon \and Juan Pablo Rossetti}

\address{Dartmouth College,  Hanover, New Hampshire, \ 03755, U.S.A.}
\email{carolyn.s.gordon@dartmouth.edu}

\address{FAMAF (ciem), Universidad Nacional de C\'ordoba, Ciudad Universitaria, 5000-C\'ordoba, Argentina.}
\email{rossetti@mate.uncor.edu}



\thanks{ The first author is partially supported by a grant from the
National Science Foundation.  The second
author is partially supported by Conicet}

\begin{document}

\maketitle

\begin{abstract}  Let $M$ be a $2m$-dimensional compact Riemannian
manifold.  We show that the spectrum of the Hodge Laplacian acting
on $m$-forms does not determine
whether the manifold has boundary, nor does it determine the
lengths of the closed geodesics.  Among the many examples are a
projective space and a hemisphere that have the same Hodge spectrum on
1-forms, and hyperbolic surfaces, mutually isospectral on 1-forms, with
different injectivity radii.  The Hodge
$m$-spectrum also does  not distinguish orbifolds from
manifolds.
\end{abstract}


\

To what extent does spectral data associated with the Hodge
Laplacian $dd^*+d^*d$
on a compact Riemannian manifold $M$ determine
the geometry and topology of $M$?  Let $\spec_p(M)$ denote the
spectrum  of the Hodge Laplacian acting on the space of $p$-forms
on $M$, with absolute boundary
conditions in case the boundary of $M$ is non-empty.  (We will
review the notion of absolute and relative boundary conditions in
Section 1.)  For each $p$, the
spectrum
$\spec_p(M)$ is known to contain considerable geometric
information.  For example, under genericity conditions, the
$p$-spectrum of a closed Riemannian manifold
$M$ determines the geodesic length spectrum of $M$.

In this article, we will focus on even-dimensional manifolds and
give a very
simple method for obtaining manifolds with the same Hodge spectrum
in the middle degree.  Through examples, we will discover that this
middle degree spectrum
contains a perhaps surprising lack of topological and geometric
information.

Among the examples of manifolds that we will construct with the
same middle degree spectrum are:

\begin{itemize}
\item  a cylinder, a M\"obius strip, and a Klein bottle;

\item  pairs of cylinders with different boundary
lengths;

\item  a hemisphere and a projective space;

\item  products $S(r)\times \P(s)$ and $\P(r)\times
S(s)$, where $\P(t)$, respectively $S(t)$, denote a projective
space, respectively sphere, of radius
$t$ (for generic $r$ and $s$, the length spectra differ);

\item non-orientable closed hyperbolic surfaces with
different injectivity radius and length spectrum;

\item  pairs of surfaces, one with an arbitrarily large
number of boundary components and the other closed.  The metrics
can be chosen to be hyperbolic.
\end{itemize}

We also obtain examples in all even dimensions exhibiting similar
properties.

These examples prove:

\begin{teor} The middle degree Hodge spectrum of an
even-dimensional Riemannian manifold $M$ does {\it not} determine:

\begin{itemize}
\item[(i)]  the volume of the boundary or even
whether $M$ has boundary;

\item[(ii)] the geodesic length spectrum or injectivity radius.  In
particular, it does  not determine the length of the shortest
closed geodesic.
\end{itemize}

\end{teor}

   The result (i) is  new.  Concerning (ii), R. Miatello and the
second author \cite{MR3} recently constructed examples of flat
manifolds, $p$-isospectral for
various choices of $p$, such that the length of the shortest closed
geodesic differed.  (Here, we say two manifolds are $p$-isospectral
if they have the same
$p$-spectrum.)

Note that all the closed manifolds that we construct are
non-orientable.  For {\it orientable} surfaces, the 1-spectrum
completely determines the 0-spectrum (which
coincides with the 2-spectrum).  Thus the behavior exhibited by our
examples cannot occur for 1-isospectral orientable surfaces.  For
{\it non-orientable} surfaces,
the 0-spectrum and 2-spectrum no longer coincide and the 1-spectrum
is the join of the two (except for the 0-eigenspace).  Our examples
show that the join of
the two spectra contains less geometric information than either one
individually.  In higher dimensions, it is perhaps surprising that
all our examples are
non-orientable.  The Hodge $*$ operator intertwines the exact and
co-exact middle degree forms and commutes with the Laplacian; thus
the ``exact'' and ``co-exact
spectra'' coincide.  Hence the $m$-spectrum of a $2m$-dimensional
{\it orientable} closed manifold, $m>1$, contains half the
information contained in the
$(m-1)$-spectrum and nothing more (except for the dimension of the
space of harmonic forms).  In contrast, for non-orientable
manifolds, the ``exact'' and
``co-exact spectra'' no longer coincide, so the middle degree
spectrum obstensibly contains the same amount of data as the other
spectra.

 One of the primary tools for recovering geometric and topological
information from the spectra is through the
small time asymptotics of the heat equation.  See, for example,
\cite{Gi} or \cite{BBG1}, \cite{BBG2}.  For closed manifolds, the trace of
the heat kernel associated with the Hodge Laplacian on $p$-forms
has an asymptotic expansion of the form
$$\zeta(t)=(4\pi t)^{-n/2}(a_0(p)+a_1(p)t+a_2(p)t^2+\dots).$$
 The coefficients $a_i(p)$ are spectral invariants and are given by
integrals, with respect to
the Riemannian measure, of expressions involving the curvature and
its covariant derivatives.  The first three heat
invariants are given as follows:
$a_0=C_0(p)\text{vol}(M)$;
$a_1(p)=C_1(p)\int_M\,\tau$, where $\tau$ is the scalar curvature
and $C_i(p)$, $i=1,2$, is a constant depending only on $p$; and
$a_2(p)$ is a linear combination, with coefficients depending only
on $p$, of $\int_M\,\tau^2$,  $\int_M\,\|Ric\|^2$ and
$\int_M\,\|R\|^2$.

For manifolds
with boundary, the expansion is instead in powers of
$t^{\frac{1}{2}}$.  The coefficient $a_{\frac{1}{2}}(p)$ is
given by $c(p)\text{vol}(\partial M)$, where $c(p)
=\binom{n-1}p-\binom{n-1}{p-1}$.

A very difficult open question is whether the heat invariants are
the only integral invariants of the spectrum,  i.e., the only
invariants which are
integrals of curvature expressions either over the manifold or over
the boundary.  This article and work of Dorothee Schueth lend
support (albeit in a
small way) towards an affirmative answer since:

\begin{itemize}
\item  The coefficient $a_{\frac{1}{2}}(p)
=c(p)\text{vol}(\partial M)$ vanishes precisely when $\dim(M)$ is
even and $p=\frac{1}{2}\dim(M)$.  Our
examples show in this case that the $p$-spectrum indeed does not
determine the volume of the boundary.

\item Dorothee Schueth
\cite{S1} showed that the three individual integral terms in the
expression for $a_2(0)$ are not spectral invariants.
\end{itemize}

A secondary theme of this paper is the spectrum of Riemannian
orbifolds.  Riemannian orbifolds are analogs of Riemannian
manifolds but with singularities.  They are locally  modelled on
quotients  of Riemannian manifolds by finite effective
isometric group actions.  The singularities  in the orbifold
correspond to the singular orbits.   One can define the notion
of Laplacian on
$p$-forms  on orbifolds.  A natural question is whether the
spectrum contains information about the singularities.  We will
see:

\begin{teor} The  middle degree Hodge spectrum cannot
distinguish Riemannian manifolds from Riemannian orbifolds with
singularities.
\end{teor}

 For example, we will see that the mutually
1-isospectral cylinder, Klein bottle and M\"obius strip are also
1-isospectral to a ``pillow'', a 2-dimensional
orbifold with four singular points, so named because of its shape.
In dimension $2m$, $m>1$, we will also construct orbifolds that are
$m$-isospectral to a projective space and to a hemisphere. To our
knowledge, the examples given here are the first examples of
orbifolds
$p$-isospectral to manifolds for some
$p$.  We do not know whether the spectrum of the Laplacian on
functions can distinguish orbifolds from manifolds, although we
will give some positive
results in
$\S 3$ and we will verify that none of the middle degree
isospectral  orbifolds and manifolds we construct are
0-isospectral.

We will also construct orbifolds, isospectral in the middle degree,
having singular
sets of arbitrarily different dimensions.

 In this paper we  have concentrated on $m$-isospectrality of
$2m$-manifolds and orbifolds.  However, for various choices of $p$,
one can also obtain
$p$-isospectral orbifolds  with singular sets of different
dimensions.  This phenomenon may be discussed in a later  paper.

We conclude these introductory remarks with a partial history of
the isospectral problem for the  Hodge Laplacian.  For a
general survey of isospectral manifolds, see \cite{Go3}.  Most
known examples of isospectral manifolds, e.g., those constructed by
the method of Sunada
\cite{Sun}, are strongly  isospectral; in particular, the manifolds
are
$p$-isospectral for all $p$.  Thus they do not reveal possible
differences in the geometric information contained in the
various $p$-spectra.  The
article
\cite{Go1} gave the first example of $0$-isospectral manifolds
which are not 1-isospectral; further examples were given in
\cite{Gt1},\cite{Gt2}.  (The $0$-spectrum, i.e., the spectrum of the
Laplacian acting
on smooth functions, is frequently referred to simply as the
spectrum of the manifold.)  Ikeda
\cite{Ik} constructed, for each positive integer
$k$,  spherical space forms  which are
$p$-isospectral for all
$p$ less than $k$ but not for $p=k$.  In the past decade, many
examples have been constructed of $0$-isospectral manifolds with
different local as well as global geometry, e.g. \cite{Sz1}, \cite{Sz2},
\cite{Go2}, \cite{Go4},\cite{GW}, \cite{GGSWW}, \cite{GSz}, \cite{S1},\cite{S2},\cite{S3}, and
\cite{Gt1};
at least in those cases in which comparisons of the higher
$p$-spectra have been carried out, these manifolds are not
$p$-isospectral for $p\geq 1$. The first
examples of manifolds which are
$p$-isospectral for some values of
$p$ but not for
$p=0$ were flat  manifolds constructed by R. Miatello and the
second author in
\cite{MR1},\cite{MR2},\cite{MR3}.   R. Gornet and J. McGowan (private communication)
recently constructed examples of spherical space forms
which are simultaneously $p$-isospectral for $p=0$ and for
various other, but not all, $p$.

The second author wishes to thank the mathematics department at
Dartmouth College for its great hospitality during the time
this paper was written.

\section {Method for constructing manifolds isospectral in the
middle degree}

\begin{definition}
Suppose $M$ is a compact Riemannian
manifold with boundary.   Let $\omega$ be a smooth $p$-form on $M$.  For
$x\in\partial M$, write
$\omega_x=\omega_x^T+\omega_x^N$, where
$\omega_x^T\in\bigwedge^p(\partial M)_x^*$ and $\omega_x^N
=\sigma\wedge\mu$ with $\sigma\in \bigwedge^{p-1}(\partial
M)_x^*$ and $\mu$ a vector in $T_x^*(M)$ normal to the cotangent
space of $\partial M$ at $x$.  The form $\omega$ is said to
satisfy {\it absolute}, respectively {\it relative}, boundary
conditions if $\omega^N=0=(d\omega)^N$, respectively
$\omega^T=0=(\delta\omega)^T$, everywhere on $\partial M$.  If
$M$ is orientable, observe that $\omega$ satisfies the relative
boundary conditions
if and only if $*\omega$ satisfies the absolute boundary
conditions, where $*$ is the Hodge duality operator.
\end{definition}

\begin{notandrem}
Let $M$ be a compact Riemannian manifold of dimension $n=2m$.

(i) If $M$ is closed let $\spm(M)$ denote the
spectrum of the Hodge Laplacian of $M$ acting on smooth $m$-forms.

(ii)  If
$M$ has boundary, we denote by  $\spm(M)$ the absolute spectrum of
$M$.  If $M$ is orientable, then since the Hodge * operator carries
$\A^m(M)$ to itself and commutes with the Hodge Laplacian, we see
from
1.1 that the absolute
$m$-spectrum coincides with the relative $m$-spectrum.  However,
this
coincidence of the spectra does not in general occur when $M$ is
non-orientable.

(iii)  Suppose that $M$ is closed and $\tau$ is an orientation
reversing
isometric involution.  If $\tau$ does not act freely, then
$\langle\tau\rangle\bs M$
is a Riemannian orbifold.  Its $m$-spectrum is defined to be the
spectrum of the
Hodge Laplacian of $M$ acting on the space of $\tau$-invariant
$m$-forms.  We denote this spectrum by $\spm(\tm)$.  (In this
article we will be concerned only with orbifolds of this form.  For
a more general discussion of
orbifolds, see for example \cite{Sc} or \cite{T}.)

If the fixed point
set of
$\tau$ is a submanifold of codimension one, then the underlying
space of the orbifold $\langle\tau\rangle\bs M$ is a Riemannian
manifold with boundary.  The
boundary corresponds to the singular set of the orbifold. (E.g., if
$M$ is a sphere and
$\tau$ is reflection across the equator, then
$\tm$ is an orbifold whose underlying space is a hemisphere.)
Since the orbifold
$m$-spectrum
$\spm(\tm)$ just defined agrees with the $m$-spectrum of the
underlying manifold with absolute boundary conditions, we will
ignore the orbifold structure and view
$\langle\tau\rangle\bs M$ as a manifold with boundary in what
follows. This point of view allows us to obtain manifolds with
boundary as quotients of closed
manifolds.
\end{notandrem}

\begin{theorem}
Let $M$ be a $(2m)$-dimensional orientable
closed Riemannian manifold.  Suppose that
$\tau$ is an orientation reversing involutive isometry of $M$.
Then
$\spm(\tm)$ consists precisely of the eigenvalues of $\spm(M)$ but
with
all multiplicities multiplied by $\frac{1}{2}$.
\end{theorem}

\begin{proof}  Let $\lambda$ be an eigenvalue in $\spm(M)$, say of
multiplicity $\ml$.
Both $\tau$ and the Hodge * operator commute with the Hodge
Laplacian and
thus leave the
$\lambda$-eigenspace $H_\lambda\subset\A^m(M)$ invariant.  Letting
$\hpl$, respectively $\hml$, denote the space of $\tau$-invariant,
respectively $\tau$-anti-invariant, forms in $\hl$, then
$\hl=\hpl\oplus\hml$.  The $\lambda$-eigenspace
in $\A^m(\tm)$ corresponds to $\hpl$.

The Hodge $*$ operator interchanges $\hpl$ and $\hml$.  To see
this, let $\Omega$ denote the Riemannian volume form of $M$.  For
$\alpha\in\A^m(M)$, we have
$\alpha\wedge *\alpha=\|\alpha\|^2\Omega$ and
$\tau^*\alpha\wedge\tau^*(*\alpha)=\|\alpha\|^2\tau^*\Omega
=-\|\alpha\|^2\Omega$,
since
$\tau$ is an orientation reversing isometry.  Thus $*$ interchanges
$\hpl$ and $\hml$.

Consequently, $\hpl$ and $\hml$ both have dimension
$\frac{1}{2}\ml$,
and the theorem follows.
\end{proof}

\begin{corollary}  Let $M_1$ and $M_2$ be
$(2m)$-dimensional orientable
closed Riemannian manifolds with $\spec_m(M_1)=\spec_m(M_2)$.
Suppose that
$\tau_1$ and
$\tau_2$ are orientation reversing involutive isometries of $M_1$
and $M_2$, respectively.  Then
$$\spm(\langle\tau_1\rangle\bs M_1)=\spm(\langle\tau_2\rangle\bs
M_2).$$
Moreover, if $N$ is a $2k$-dimensional closed Riemannian manifold,
then
$$\spec_{m+k}((\langle\tau_1\rangle\bs M_1)\times
N)=\spec_{m+k}((\langle\tau_2\rangle\bs M_2)\times N).$$
\end{corollary}

For the second statement, observe that $\tau_i$ extends to an
involutive orientation reversing isometry of $M_i\times N$.

\begin{weakremark}
The conclusion fails to hold if we drop the
hypothesis that the involutions be orientation reversing.  For
example, Let $M_1=M_2=(\Z\times
2\Z)\bs\R^2$, and let $\tau_1$ and $\tau_2$ be the translations
$(x,y)\mapsto(x+\frac{1}{2},y)$ and $(x,y)\mapsto (x,y+1)$,
respectively. Then the quotient tori are
not 1-isospectral.
\end{weakremark}

In our applications of Corollary 1.4 below, we will take $M_1
=M_2$.

A useful special case of Corollary 1.4 is the following:

\begin{corollary}  Let $M$ and $M'$ be manifolds of
dimension $k$ and $k'$, respectively, with $k+k'$ even, say $k+k'=2m$.  
Suppose $\phi$ and $\phi'$ are
orientation reversing involutive  isometries of $M$ and $M'$
respectively.  Then $(\langle\phi\rangle\backslash M)\times M'$ is
$m$-isospectral to $M\times
(\langle\phi'\rangle\backslash M')$.
\end{corollary}

Given a Riemannian manifold $M$ and a non-negative integer $k\leq\dim(M)$, let $\{\lambda_j\}_{j=1}^\infty$ be 
the spectrum of the Laplacian acting on $k$-forms on $M$.  The complexified heat trace is 
given by $Z(z)=\Sigma_{n=0}^\infty\,\exp(-\frac{\lambda_n}{z})$ for $z\in\C$.  
Recall that Y. Colin de Verdi\`ere \cite{C} obtained an expansion of the complexifed 
heat trace as a sum of oscillating terms whose periods are related to the lengths of 
the closed geodesics in $M$.  The authors are grateful to the referee for pointing out 
the following additional corollary of Theorem~1.3.

\begin{corollary}  Suppose that $N$ is a nonorientable closed Riemannian manifold of
dimension $2m$.  Let $\gamma$ be an isolated, nondegenerate simple closed geodesic such
that the holonomy about $\gamma$ is orientation reversing.  Then $\gamma$ does not contribute
to the complexified heat trace associated with the Laplacian acting on $m$-forms.
\end{corollary}

\begin{proof}  Let $M$ be the orientation covering of $N$ with the lifted Riemannian metric.
By Theorem 1.3, the complexified heat traces of $M$ and $N$ associated with the Laplacian
on $m$-forms coincide except for a factor of $\frac{1}{2}$.  Since the lift of $\gamma$ to $M$
is not closed, it does not contribute to the complexified heat trace of $M$ and hence
$\gamma$ does not contribute to that of $N$.
\end{proof}

\begin{remark} In the setting of Corolllary 1.6 with $m=1$, it was already known that the contribution 
of $\gamma$ to the complexified heat trace on $1$-forms must be zero up to principal order since the 
principal order coefficient is the trace of the holonomy map about the geodesic.
\end{remark}

\section{Examples of manifolds isospectral in the middle degree}

\begin{remark} Most, though not all, of the examples we
construct below will be surfaces.  However, by taking  products of
the manifolds in these examples
with a closed manifold
$N$ and applying the second statement of Corollary 1.4, one obtains
examples in arbitrary even dimensions.
\end{remark}

We begin with manifolds of positive curvature.

\begin{nada}{\bf  Hemispheres and projective spaces.}  Let $M$
be the $(2m)$-sphere, $m\geq 1$, let $\tau$
denote the antipodal map and let $\sigma$ denote reflection about
an equatorial sphere.  Then Corollary 1.4 implies that the projective
space $\tm$ is $m$-isospectral to the hemisphere $\sm$. (I.e., the
$m$-spectrum of the projective space coincides with the absolute
$m$-spectrum of the hemisphere.)
\end{nada}

\begin{nada}{\bf  Half ellipsoids with different boundary volume.}
Let $M$ be a $2m$-dimensional ellipsoid (not round) and let
$\sigma$ and $\tau$
denote the reflections across two different hyperplanes of
symmetry.  Then
$\tm$ and
$\sm$ are $m$-isospectral manifolds with boundary, but their
boundaries
have different volume.
\end{nada}

\begin{nada}{\bf Products of spheres and projective spaces with
different length spectra and injectivity radius}  Let $S^k(r)$ and
$\P^k(r)$ denote the sphere
and projective space, respectively, of  radius $r$ and dimension
$k$.  Then by Corollary 1.5, for $m,n\geq 1$ and for any $r,s>0$,
the manifolds $S^{2m}(r)\times
\P^{2n}(s)$ and
$\P^{2m}(r)\times S^{2n}(s)$ have the same $(m+n)$-spectrum.
However, for generic choices of $r$ and $s$, the lengths of the
shortest closed geodesics in the two
manifolds differ.

These closed manifolds are also  isospectral to the product of a
hemisphere and a sphere, as can be seen by taking an
involution of $S^{2m}(r)\times S^{2n}(s)$
given by the equatorial reflection in one factor
\end{nada}

We next consider flat manifolds.

\begin{nada}{\bf Cylinders, Klein bottles, and M\" obius strips.}
The following, with appropriately chosen sizes, are mutually
1-isospectral:
\begin{itemize}
\item[(i)]   a Klein bottle;
\item[(ii)]  a cylinder, say of perimeter $p$;
\item[(iii)] a M\" obius strip of perimeter $\frac{p}{\sqrt{2}}$.
\end{itemize}

Indeed each of these is a quotient of $T=\Z^2\bs \rt$ by an
involution $\tau$ arising from an involution $\tilde{\tau}$ of
$\rt$ given as follows:
For the Klein bottle, $\tilde{\tau}$ is the composition of
translation by $(\frac{1}{2},0)$ with
reflection $R$ about the $x$-axis.
For the cylinder,  $\tilde{\tau}$ is the reflection $R$ just
defined.
For the M\"obius strip, $\tilde{\tau}$ is reflection across the
line $x=y$.

Thus by Corollary 1.4, the Klein bottle, the cylinder and the
M\"obius strip all have the same 1-form
spectrum.
\end{nada}

\begin{nada}{\bf Cylinders of different shapes.}
Let $C_{(h,w)}$ denote the flat cylinder of height $h$ and
circumference $w$, hence of perimeter $2w$.  Then for
all $a,b>0$, we have
$$\spec_1(C_{(a,\frac{b}{2}})=\spec_1(C_{(b,\frac{a}{2})})$$
since both cylinders are quotients by a reflection of a rectangular
torus of height $a$ and
width $b$.
Thus cylinders of arbitrarily different perimeters can have the
same 1-form spectrum.
\end{nada}

\begin{remarks}
(i) One similarly obtains pairs of non-isometric Klein bottles with
the same 1-spectra and with different lengths of closed geodesics.

(ii) By taking $M$ to be a rhomboidal torus and letting $\tau_1$
and $\tau_2$ be reflections across the ``diagonals'' one similarly
obtains
pairs of 1-isospectral M\"obius strips of arbitrarily different
perimeters.

(iii)  By taking products of the various 1-isospectral flat
manifolds in 2.5, 2.6 and in these remarks with a torus of
dimension $2m-2$ and recalling Remark
2.1, we obtain pairs of
$m$-isospectral $2m$-dimensional flat manifolds with the same
boundary behavior as the surfaces in these examples.
\end{remarks}

Finally we consider hyperbolic surfaces.

\begin{theorem}  In each genus $g>1$, there exist
non-orientable hyperbolic closed surfaces $S_1$, $S_2$, $S_3$, and
$S_4$ of genus $g$ such that:

\item{(a)}  All four surfaces have the same Hodge spectrum on
$1$-forms;

\item{(b)}  They have different length spectra (different in the
strong sense that different lengths occur, not just different
multiplicities), and $S_4$ has a
smaller injectivity radius than
$S_1$,
$S_2$, and
$S_3$.

\item{(c)}  The Laplacians acting on functions on these surfaces
are not isospectral.

 These surfaces are also 1-isospectral to each of four hyperbolic
surfaces with boundary having 2, 4, $2g-2$, and $2g$ boundary
components, respectively.
\end{theorem}

We remark that (b) implies (c), since the 0-spectrum determines the
weak length spectrum (i.e., the spectrum of lengths of closed
geodesics, ignoring
multiplicities) in the case of hyperbolic manifolds.  (Aside:  For
Riemann surfaces, it is a classical result of Huber that the
0-spectrum and the strong length
spectrum, i.e., the spectrum of lengths, counting multiplicities,
determine each other.  This result has been extended to
non-orientable hyperbolic surfaces by
Peter Doyle and the second author in an article in preparation.)

\begin{figure}[!htb]
\input{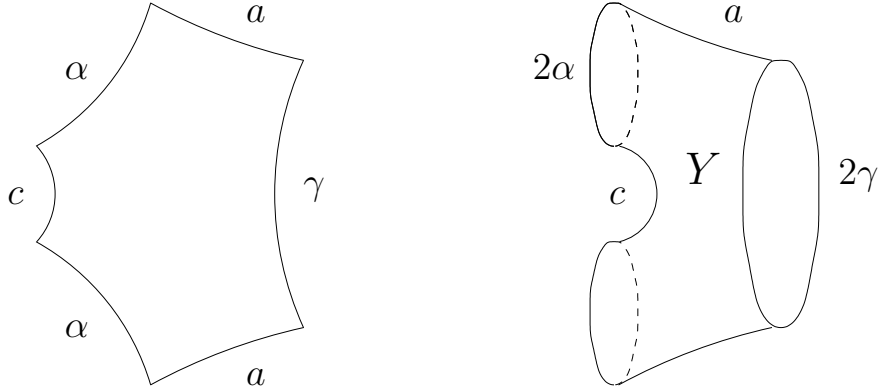}
\caption{Right-angled hexagon and the associated pair of pants.}
\label{fig1}
\end{figure}

\begin{proof}  We construct an orientable surface $S$ as follows:
Let $Y$ be a pair of hyperbolic  pants, as shown in Figure 1; the
boundary geodesics of the pant legs have the same length, while the
waist may have a different length.  The pants are formed by gluing
together two identical right-angled geodesic hexagons along three
sides.  The hexagon is also pictured in Figure 1.  Choose a positive
integer
$t$ and glue together
$4t$ isometric copies
$Y_1,\dots,Y_{4t}$ of
$Y$  to obtain a Riemann surface $S$ of genus $2t+1$. To describe
the symmetries, we will visualize $S$
(as shown in Figure 2 in case $t=1$) as obtained from a
surface $N$ with boundary in  $\R^3$ with appropriate
identifications of boundary edges.
The three  symmetries $\tau_H$,
$\tau_{V1}$, and $\tau_P$  of $N$ given by reflection across the
$xy$-plane, the $xz$-plane and the $yz$-plane, respectively, define
commuting involutive
isometries of
$S$.  (We are using the indices $H$,
$V$, and $P$ here to indicate reflection across planes that are
horizontal, vertical, or the plane  of the paper.)  There are
additional isometric involutions of
$S$ (not of $N$) given by reflections across vertical planes
passing through the waists of two of the pants.  Choose such a
symmetry and denote it by
$\tau_{V2}$.  E.g., in case $t=1$, a choice of $\tau_{V2}$
 interchanges pants $Y_1$ with $Y_2$ and $Y_3$ with $Y_4$ in Figure 2.
Finally, let $\rho$ be the orientation preserving isometric
involution (rotation)
$\rho$ sending the pair of pants $Y_i$ to the pants $Y_{i+2t}$ mod
$4t$, for each $i$.  Note that $\th$ and $\tb$ commute
with
$\rho$.  Thus the isometries $\tone :=\th\circ\tvo\circ\tb$,
$\ttwo:=\tb\circ\rho$, $\tthree:=\th\circ\rho$, and $\tfour :
=\th\circ\tvt\circ\tb$ are involutive, orientation-reversing,
fixed-point-free isometries of $S$.  We set
$S_i=\langle\tau_i\rangle\backslash S$, for $i=1,2,3, 4$.  The surfaces
$S_i$  are depicted in Figure 3 in case $t=1$.  These
non-orientable surfaces have genus $g=t+1$, and,
topologically, they are spheres
with $t-1$ handles and two cross handles, or equivalently, spheres
with $t+1$
cross handles.

\begin{figure}[!htb]
\input{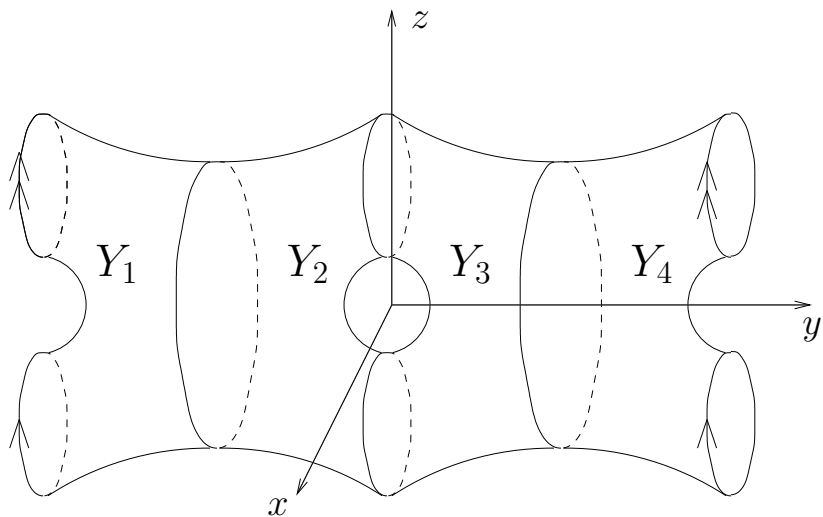}
\caption{Four congruent pants glued to form $S$.}
\label{fig2}
\end{figure}

For generic choices of the hyperbolic pants $Y$ used 
to build the surfaces above, the four surfaces will 
have different geodesic length spectra.  For 
concreteness, we will make a specific
choice of $Y$ to guarantee that one of the surfaces has
a strictly smaller injectivity radius than the others. 
There exists a unique right-angled hexagon as in Figure 1 for each 
given choice of $\alpha$ and $\gamma$.  (See Buser \cite{B},
Theorem 2.4.2.)   We choose $\alpha$ and 
$\gamma$ to be less than $\text{arcsinh}(1)$ and to satisfy 
$\gamma <2\alpha$.  Since the waist of the resulting pants
has length $2\gamma$ and the boundary geodesics of the 
legs have length $2\alpha$, the resulting surface $S$
has $2t$ closed geodesics of length $2\gamma$ and $4t$
of length $2\alpha$.  By Theorem 4.1.6 of \cite{B}, a surface
of genus $g$ can have at most  $3g-3$ simple closed 
geodesics of length less than or equal to $2 \text{arcsinh}(1)$.
Thus the geodesics in $S$ corresponding to the boundary 
geodesics of the pants are the only simple closed geodesics
satisfying this bound on their lengths.  
  
Now for $i=1,2,3,4$, consider the shortest ``new''
closed geodesic in the surface $S_i$, i.e., the shortest
closed geodesic that does not lift to a closed geodesic of the same length
in $S$.  These geodescis are depicted in Figure 3 in case $t=1$.  
(The relative lengths of the geodesics in the figures 
are distorted.)   The surface $S_4$ contains a geodesic 
of length $\gamma$.  This geodesic is shorter than any 
closed geodesic in the orther surfaces.  Thus $S_4$ has
smaller injectivity radius than the others.
Thus statement (b) of the Theorem is satisfied.

\begin{figure}[!htb]
\input{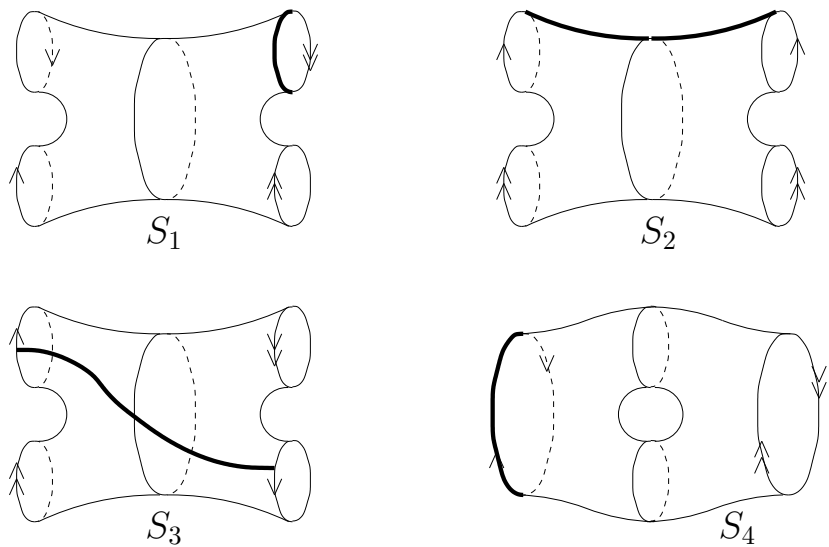}
\caption{The non-orientable surfaces $S_i$; the shortest ``new'' geodesics are thicker.}
\label{fig3}
\end{figure}

For the final statement of the theorem, we use Corollary 1.4 to
observe that the surfaces $S_i$ are also 1-isospectral to the four
surfaces $\langle\tau_P\rangle\bs S$,
$\langle\tau_H\rangle\bs S$, $\langle\tau_{V1}\rangle\bs S$, and
$\langle\tau_{V2}\rangle\bs S$ with  boundary.  These surfaces have $2g$,
$2g-2$, 4 and 2 boundary  components, respectively.
\end{proof}

\begin{remark}  Letting $S$ be the surface of genus $2n$,
$n\geq 1$, pictured in Figure 4, one may take the quotient of $S$
by each of the three visual
symmetries (reflection across a vertical plane, horizontal plane,
and plane of the paper) to obtain 1-isospectral surfaces with
boundary.  The first has only one
boundary component, while the others have $2n+1$.
\end{remark}

\begin{figure}[!htb]
\input{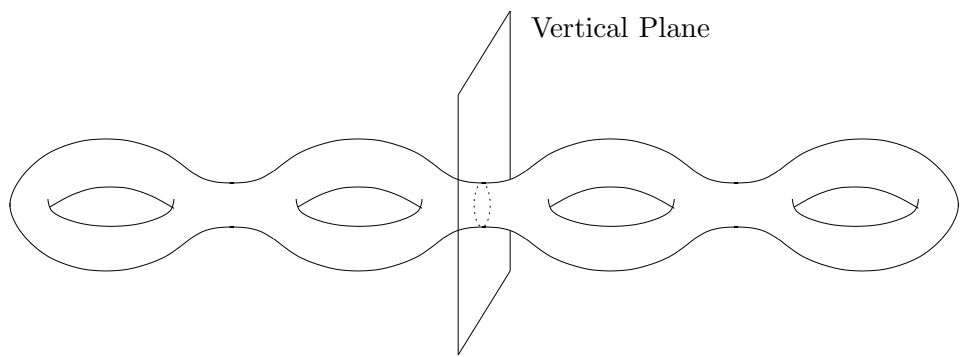}
\caption{A Riemann surface and a reflection plane.}
\label{fig4}
\end{figure}

\begin{nada}{\bf Surfaces with boundaries of arbitrarily different
shapes.}  This final example involves manifolds of mixed curvature.
Let $M$ be
a 2-sphere with a non-standard metric, symmetric with respect to
reflection across two orthogonal planes.  The quotients of $M$ by
these isometric involutions are
1-isospectral manifolds with boundaries.  By choosing the metrics
on $M$ appropriately, one can more or less arbitrarily prescribe
independently the geodesic
curvatures of the boundary curves.
\end{nada}

\section{Examples of orbifolds and manifolds isospectral in the
middle degree}

Does the spectrum distinguish orbifolds with singularities from
smooth manifolds?
We will see below that the middle degree Hodge spectrum cannot
detect the presence of singularities.  We do not know whether the
0-spectrum always detects
singularities.  However, we will show in Propositions 3.4 and 3.5
that it does in many situations, including all the examples we will
give of orbifolds and
manifolds which are isospectral in the middle degree.

\begin{theorem}  The cylinder, Klein bottle and M\"obius
strip of Example 2.5 are also 1-isospectral to a four pillow $\O$.
\end{theorem}

\begin{proof} The four pillow is the quotient of $T=\Z^2\bs \rt$ by
the involution $\rho$ induced from the map $-Id:\R^2\to\R^2$.  It
is a pillow-shaped orbifold with
isolated singularities at the four corners.  Since
$\rho$ is not orientation reversing,  we cannot apply Corollary
1.4.  Instead we will compare the spectra directly.

Every 1-form on $T$ may be written in the form $f\,dx +g\,dy$,
where
$dx$ and
$dy$ are the forms induced on
$T$  by the standard forms on
$\rt$.  Since
$$\Lap(f\,dx+g\,dy)=\Lap(f)\,dx+\Lap(g)\,dy,$$
$\spec_1(T)$ consists of two copies of $\spec_0(T)$.   Hence each
of the manifolds in Example 2.5 has 1-spectrum equal to the
0-spectrum of $T$.

Since both $dx$ and $dy$ are
$\rho$-anti-invariant, the space of $\rho$-invariant one-forms on
$T$ is given by
$$\{f\,dx +g\,dy:\tau^*f=-f, \tau^*g=-g\}. $$
Thus $\spec_1(\O)$ is formed by two copies of the spectrum of the
Laplacian of $T$ acting on the space
of
$\rho$-anti-invariant smooth functions.  By Fourier decomposition
on the torus, one sees that the Laplacian of $T$ restricted to the
$\rho$-anti-invariant smooth
functions has the same spectrum as the Laplacian restricted to the
$\rho$-invariant smooth functions.  Consequently
$\spec_1(\O)=\spec_0(T)$ and the theorem follows.
\end{proof}

\begin{theorem} (i)  For each  positive integer $m>1$,
there exists a collection of $m-1$ distinct mutually
$m$-isospectral
$2m$-dimensional flat orbifolds having
singular sets of dimension $1, 3, 5, \dots, 2m-3$, respectively.
Moreover, these orbifolds are also m-isospectral to the direct
product of a cubical
$(2m-2)$-torus with each of the manifolds of Example 2.5.

(ii)  For each  positive integer $m>1$, there exists a collection
of $m-1$ distinct mutually $m$-isospectral
$2m$-dimensional spherical orbifolds having
singular sets of dimension $ 1, 3, \dots, 2m-3$, respectively.
Moreover, these orbifolds are also m-isospectral to a
$2m$-dimensional projective space and to a
$2m$-dimensional hemisphere.
\end{theorem}

\begin{weakremark}  One can also incorporate into the family in
Theorem 3.2(i) of mutually $m$-isospectral orbifolds and manifolds
a $2m$-dimensional analogue
of the  pillow, obtained as the quotient of $T^{2m}$ under the
action of minus the identity. This orbifold has only isolated
singularities; i.e., the dimension of its singular set is zero. The
proof of isospectrality is similar to that of Theorem~3.1.
\end{weakremark}

\begin{proof} (i) For $k$  odd with $1\leq k \le 2m-3$, let $\tau_k$
be the orthogonal involution  of $\R^{2m}$ given by the $2m\times
2m$ diagonal matrix
$\text{diag}
(\underbrace{-1,\dots,-1}_{2m-k},\underbrace{1,1,\dots,1}_{k})$.
Then $\tau_k$
induces an orientation reversing involution,  which we also denote
by $\tau_k$, of the
$2m$-dimensional cubical flat torus $T$.  Let $\O_k$ be the
quotient of $T$ by $\tau_k$.  By Corollary 1.4, these orbifolds are
mutually $m$-isospectral.  Since
the direct product of a cubical
$(2m-2)$-torus with any of the manifolds or orbifolds in Theorem
3.1 may also be viewed as a quotient of $T$ by an involution, the
final statement of the theorem
follows.  (Aside:  One could also take $k=2m-1$ to obtain an
orbifold with singular set of dimension $2m-1$.  However, under the
identification described in
Remark 1.2(iii), this orbifold is identified with the product of
the cyliner in Theorem 3.1 with the cubical torus and thus is
redundant.)

(ii) We now let $\tau_k$, for $k$  even with $2\le k \le 2m-2$, be
the orthogonal involution  of $\R^{2m+1}$ given by the
$(2m+1)\times (2m+1)$ diagonal matrix
$\text{diag}
(\underbrace{-1,\dots,-1}_{2m+1-k},\underbrace{1,1,\dots,1}_{k})$.
The quotient of the $2m$-sphere by $\tau_k$ is
an orbifold with singular set of dimension $k-1$.  Again by
Corollary 1.4, these orbifolds are mutually $m$-isospectral and are
also $m$-isospectral to the
projective space and hemisphere of Example 2.2.  (Aside:  If we
take $k=2m$, we obtain the hemisphere under the identification in
1.2.)
\end{proof}

In \cite{MR1},\cite{MR2}, R. Miatello and the second author computed all the
various $p$-spectra of flat manifolds.  Their method extends to
flat orbifolds and can be used
to give an alternative proof of Theorem 3.2(i) as well as of the
earlier examples of flat manifolds and/or orbifolds.

We next show that the orbifolds and  manifolds in Theorems 3.1 and
3.2 are not 0-isospectral.   The orbifolds in these examples belong
to the class of so-called
``good'' Riemannian orbifolds; that is, they are quotients of
Riemannian manifolds by discrete, effective, properly discontinuous
isometric group actions.

\begin{lemma}\cite{D1},\cite{D2}  Let $\O$ be a good Riemannian
orbifold, and let $\lambda_1, \lambda_2,\dots$ be the 0-spectrum of
the Laplacian.  Then there is an
asymptotic expansion as $t\downarrow 0$ of the form
$$\Sigma_{i=1}^\infty\,e^{-\lambda_it}\sim (4\pi
t)^{-\frac{n}{2}}\,\Sigma_{k
=0}^\infty\,a_kt^k\,+\,\Sigma_S\,B_S(t),$$
where $S$ varies over the strata of the singular set and where
$B_S(t)$ is of the form $$B_S(t)=(4\pi
t)^{-\frac{\dim(S)}{2}}\,\Sigma_{k=0}^\infty\,b_{k,S}t^k$$
with
$b_{0,S}\neq 0$.  The
coefficients
$a_k$ in the first part of the expansion are given by the same
curvature integrals as in the heat expansion for manifolds; in
particular,
$a_0=vol(\O)$.
\end{lemma}

\begin{weakremarks} (i) There is an apparent typographical error in
the statement of this result as Theorem 4.8 in the article
\cite{D2} in that the power $(4\pi
t)^{-\frac{\dim(S)}{2}}$ is missing.  However, the proof makes the
expression clear.

(ii)  S. Greenwald, D. Webb, S. Zhu and the first author recently
generalized Lemma 3.3 to arbitrary Riemannian orbifolds.  An
article is in preparation.  With
this generalization, the hypothesis that $\O$ be good in the first
statement of Proposition 3.4 below can be dropped.
\end{weakremarks}

\begin{proposition}  Let $\O$ be a good Riemannian orbifold
with singularities.  Then:
\item{(i)}  If $\O$ is even-dimensional (respectively,
odd-dimensional) and some strata of the singular set is
odd-dimensional (respectively, even-dimensional),
then $\O$ cannot be 0-isospectral to a Riemannian  manifold.

\item{(ii)}  If $N$ is a manifold such that $\O$ and $N$ have a
common Riemannian cover, then $M$ and $\O$ are not 0-isospectral.
\end{proposition}

\begin{proof}   (i)  In the two cases, the fact that $\O$ is an
orbifold can be gleaned from the presence of half-integer powers,
respectively integer powers,
of
$t$ in the asymptotic expansion in Lemma 3.3.

(ii)  Suppose $N$ and $\O$ have a common covering.  By Lemma 3.3,
the spectrum of $\O$ determines its volume.  The first part of the
heat expansion in Lemma 3.3 (involving the
$a_k$) depends only on the volume of the orbifold and the curvature
of the covering manifold; thus it must be identical for $N$ and
$\O$.  The second part of the
expansion vanishes for $N$ but not for $\O$.  Hence the trace of
the heat kernels of $N$ and $\O$ have different asymptotic
expansions, so $N$ and $\O$ can't be
0-isospectral.
\end{proof}

\begin{proposition}  None of the manifolds and orbifolds in
Theorems 3.1 and 3.2 are 0-isospectral.
\end{proposition}

\begin{proof} By
Proposition 3.4(ii), the closed manifolds in these theorems are not
0-isospectral to the orbifolds.  As discussed in Remark 1.2, the
manifolds with boundary in
these theorems are the underlying spaces of (closed) orbifolds; the
boundaries of the manifolds form the singular sets of the
orbifolds.  Moreover the Neumann
spectrum of each such manifold coincides with the 0-spectrum of the
associated orbifold.  Thus we are reduced to comparing the
0-spectra of various orbifolds, each
of which has a singular set consisting of a single strata.
Moreover, the orbifolds in each collection have the same constant
curvature and have a common covering
as orbifolds.  In particular, the first part of the heat expansion
in Lemma 3.3 is identical for the various orbifolds.  However,
since their singular strata have
different dimensions, the second part of the expansion differs.
Thus they cannot be 0-isospectral.
\end{proof}

\end{document}